\documentclass[a4paper, 12pt]{amsproc}

\usepackage{fullpage}
\usepackage{amsmath, amsthm, amssymb, mathtools, mathrsfs}

\usepackage{hyperref}

\usepackage{chessfss}
\usepackage{tikz}
\usetikzlibrary{intersections, calc, arrows.meta}

\usepackage{graphicx}
\usepackage{here}
\usepackage{time}

\usepackage{xcolor}
\usepackage[capitalize,nameinlink,noabbrev,nosort]{cleveref}
\hypersetup{
	colorlinks=true,       
	linkcolor=brown,          
	citecolor=brown,        
	filecolor=brown,      
	urlcolor=brown,           
}

\makeatletter
\@namedef{subjclassname@2020}{%
  \textup{2020} Mathematics Subject Classification}
\makeatother

\newtheorem{theoremcounter}{Theorem Counter}[section]

\theoremstyle{definition}

\theoremstyle{plain}

\newtheorem{theorem}[theoremcounter]{Theorem}

\numberwithin{equation}{section}

\def\cR{\mathscr R}

\begin{document}

\title[]{Multiple zeta values and Euler's reflection formula for the gamma function} 

\author{Karin Ikeda} 
\address{Joint Graduate School of Mathematics for Innovation, Kyushu University,
Motooka 744, Nishi-ku, Fukuoka 819-0395, Japan}
\email{ikeda.karin.236@s.kyushu-u.ac.jp}

\author{Mika Sakata}
\address{School of Physical Education, Osaka University of Health and Sport Sciences,
1-1 Asashirodai, Kumatori-cho, Sennan-gun, Osaka 590-0496, Japan}
\email{m.sakata@ouhs.ac.jp} 

\thanks{The first author was supported by WISE program (MEXT) at Kyushu University. 
}



\maketitle

\begin{abstract}
	In this paper, we give a purely algebraic proof of an identity coming directly from Euler's reflection formula for the gamma function. Our proof uses Hoffman's harmonic algebra and some binomial identities.
\end{abstract}

\section{Introduction}
For given positive integers $k_1,\ldots,k_r$ with $k_r\geq2$, the multiple zeta value (MZV) $\zeta(k_1,\ldots,k_r)$ is defined by
$$
\zeta(k_1,\ldots,k_r)=\sum_{0<m_1<\cdots<m_r}\frac{1}{{m_1}^{k_1}\cdots{m_r}^{k_r}}.
$$

Recall the Taylor series for the (log of) gamma function
$$
\Gamma(1+x)=\exp\left(-\gamma x+\sum_{n=2}^{\infty}\frac{(-1)^{n}}{n}\zeta(n)x^n\right) (\gamma:\text{Euler's constant}).
$$
This and the well-known reflection formula
$$
\frac{1}{\Gamma(1+x)\Gamma(1-x)}=\frac{\sin{\pi x}}{\pi x}
$$
together with the expansion (which is a direct consequence of the infinite product of $\sin x$)
$$
\frac{\sin{\pi x}}{\pi x}=1+\sum_{n=1}^{\infty}(-1)^n\zeta(\underbrace{2,\ldots,2}_{n})x^{2n}
$$
give
$$
\exp\left(\sum_{n=2}^{\infty}\frac{(-1)^{n-1}}{n}\zeta(n)x^n\right)\exp\left(-\sum_{m=2}^{\infty}\frac{\zeta(m)}{m}x^m\right)=1+\sum_{n=1}^{\infty}(-1)^n\zeta(\underbrace{2,\ldots,2}_{n})x^{2n}.
$$
The aim of this short note is to give a purely algebraic proof of this identity. More precisely, we prove the identity (\ref{main}) (Section 2) in the framework of ``Hoffman's algebra''.
\section{Algebraic setup and the main theorem}
We recall the algebraic setup of MZVs introduced by Hoffman \cite{Hof97}. We are dealing directly with indices, not with non-commutative polynomials as in \cite{Hof97}. Let $\cR$ be the $\mathbb{Q}$-vector space
$$
\cR=\bigoplus_{r=0}^{\infty}\mathbb{Q}[\mathbb{N}^r]
$$
spanned by finite $\mathbb{Q}$-linear combinations of symbols $[\textbf{k}]=[k_1,\ldots,k_r]$ with $\textbf{k}=(k_1,\ldots,k_r)\in\mathbb{N}^r$ for some $r$. We understand $\mathbb{Q}[\mathbb{N}^0]=\mathbb{Q}[\varnothing]$ for $r=0$. Further let $\cR^0$ denote the subspace of $\cR$ spanned by the $\textit{admissible}$ symbols, i.e., by $[\varnothing]$ and the symbols $[k_1,\ldots,k_r]$ with $k_r\geq2$. On $\cR$, we consider the $\mathbb{Q}$-bilinear harmonic (stuffle) product $\ast$ which is defined  inductively as:
\begin{enumerate}
\item[(a)] for any index $\textbf{k}, [\varnothing]\ast[\textbf{k}]=[\textbf{k}]\ast[\varnothing]=[\textbf{k}]$;
\item[(b)] for any indices $\textbf{k}=(k_1,\ldots,k_r)$ and $\textbf{l}=(l_1,\ldots,l_s)$ with $r,s\geq1$,
\end{enumerate}
$$
[\textbf{k}]\ast[\textbf{l}]=[[\textbf{k}\_]\ast[\textbf{l}],k_r]+[[\textbf{k}]\ast[\textbf{l}\_],l_s]+[[\textbf{k}\_]\ast[\textbf{l}\_],k_r+l_s],
$$
where $\textbf{k}\_=(k_1,\ldots,k_{r-1}),\textbf{l}\_=(l_1,\ldots,l_{s-1})$.

Hoffman proved that $\cR_{\ast}:=(\cR,\ast)$ is a commutative and associative $\mathbb{Q}$-algebra and that $\cR_{\ast}^{0}$ is a subalgebra of $\cR_{\ast}$ \cite{Hof97}. Moreover, he proved that the evaluation map $\zeta : \cR_{\ast}^0\ni[k_1,\ldots,k_r]\mapsto\zeta(k_1,\ldots,k_r)\in\mathbb{R}$, being extended $\mathbb{Q}$-linearly, is an algebra homomorphism from $\cR_{\ast}^0$ to $\mathbb{R}$.

Our main theorem is the following.
\begin{theorem}
The equality 
\begin{align}\label{main}
&\exp_{\ast}\left(\sum_{n=2}^{\infty}\frac{(-1)^{n-1}}{n}[n]x^{n}\right)\exp_{\ast}\left(-\sum_{m=2}^{\infty}\frac{[m]}{m}x^{m}\right)=1+\sum_{n=1}^{\infty}(-1)^n[\underbrace{2,\ldots,2}_{n}]x^{2n}
\end{align}
holds in be the ring of formal power series $\cR_{\ast}[[x]]$ over $\cR_{\ast}$, where $\exp_{\ast}$ is the exponential $\exp_{\ast}(f)=\sum_{n=0}^{\infty}\frac{f^n}{n!}$ with $f^n$ being the power in the ring  $\cR_{\ast}[[x]]$.
\end{theorem}
\section{Proof} \label{s2}
Recall the equality
$$
\exp\left(\sum_{n=2}^{\infty} \frac{(-1)^{n-1}}{n}\zeta(n) x^n \right)=1+\sum_{k=2}^{\infty}(-1)^kS_{\mathbb{R}}(k)x^k,
$$
where
$$
S_{\mathbb{R}}(k)=\sum_{\substack{k_1+\cdots+k_r=k\\r\geq1,^\forall k_i\geq2}}(-1)^r\prod_{j=1}^{r}\frac{(k_j-1)}{k_j!}\zeta(k_1,\cdots,k_r).
$$
As indicated in \cite{AK}, this can be proved easily by using the Weierstrass infinite product of $\Gamma(1+x)$.

The second named author proved the algebraic counterpart of this identity. 
\begin{theorem}[\cite{Sakata}]\label{second}
The equality
\begin{align}
\exp_{\ast}\left(\sum_{n=2}^{\infty}\frac{(-1)^{n-1}}{n}[n]x^{n}\right)
=\ 1+\sum_{k=2}^{\infty}(-1)^{k}S(k)x^{k} 
\end{align}
holds in $\cR_{\ast}[[x]]$, where $S(k)=\sum_{\substack{k_{1}+\cdots +k_{r}=k\\r\geq1, ^{\forall}k_{i}\geq 2}}(-1)^{r}\prod_{j=1}^{r}\frac{(k_{j}-1)}{k_{j}!}
[k_{1},\ldots,k_{r}]$.
\end {theorem}

Using this, our task is to show the identity
\begin{align}
\left(\sum_{k=0}^{\infty}(-1)^kS(k)x^k\right)\ast\left(\sum_{l=0}^{\infty}S(l)x^l\right)=1+\sum_{n=1}^{\infty}(-1)^n[\underbrace{2,\ldots,2}_{n}]x^{2n}
\end{align}
in $\cR_{\ast}[[x]]$, where $S(0)=1, S(1)=0$.

If $n$ is an odd number, then $\sum_{k=0}^{n}(-1)^kS(k)\ast S(n-k)$ is equal to $0$ by the obvious symmetry. Hence we have
\begin{align}
\text{(LHS)}
&=\sum_{n=0}^{\infty}\sum_{k=0}^{2n}(-1)^kS(k)\ast S(2n-k)x^{2n}.
\end{align}
The proof will then be done if we show 
\begin{align}\label{eq2}
\sum_{k=0}^{2n}(-1)^kS(k)\ast S(2n-k)=(-1)^n[\underbrace{2,\ldots,2}_{n}]
\end{align}
for $n\geq1$. 

Let $(k_1,\ldots,k_l)$ be any fixed index of weight $2n$. We now compute the coefficient of $[k_1,\ldots,k_l]$ on the left-hand side of (\ref{eq2}). Let $h$ be the number of components in $[k_1,\ldots,k_l]$ which are greater than $3$. We first consider the case $h>0$. For $i$ with $0\leq i\leq h$, suppose that exactly $i$ components in $[k_1,\ldots,k_l]$, say $k_1,\cdots,k_i$, result from the harmonic product of the form (for instance)
$$
[m_1,\ldots,m_i,k_{i+1},\ldots,k_t]\ast[k_1-m_1,\ldots,k_i-m_i,k_{t+1},\ldots,k_l],
$$
where $2\leq m_j\leq k_j-2$ for $1\leq j\leq i$ (that is, $k_j=m_j+(k_j-m_j)$ for $1\leq j\leq i$) and $i\leq t \leq l$. 

Then the coefficient of $[k_1,\ldots,k_l]$ coming as a result of this product on the left-hand side of (\ref{eq2}) is, noting that the factor $[m_1,\ldots,m_i,k_{i+1},\ldots,k_t]$ may come either from $S(k)$ or $S(2n-k)$,
\begin{align*}
&(-1)^{m_1+\cdots+m_i+k_{i+1}+\cdots+k_t}(-1)^{t}\prod_{j=1}^{i}\frac{m_j-1}{m_j!}\prod_{d=i+1}^{t}\frac{k_d-1}{k_d!}\\
&\qquad\times(-1)^{i+l-t}\prod_{j=1}^{i}\frac{k_j-m_j-1}{(k_j-m_j)!}\prod_{d=t+1}^{l}\frac{k_d-1}{k_d!}\\
&+(-1)^{k_1-m_1+\cdots +k_i-m_i+k_{t+1}+\cdots+k_l}(-1)^{i+l-t}\prod_{j=1}^{i}\frac{k_j-m_j-1}{(k_j-m_j)!}\prod_{d=t+1}^{l}\frac{k_d-1}{k_d!}\\
&\qquad\times(-1)^{t}\prod_{j=1}^{i}\frac{m_j-1}{m_j!}\prod_{d=i+1}^{t}\frac{k_d-1}{k_d!}.
\end{align*}
Summing up over all the possibilities of $m_j$, we have 
\begin{align*}
&(-1)^{l+i}\sum_{m_1=2}^{k_1-2}\cdots\sum_{m_i=2}^{k_i-2}\\
&\left\{(-1)^{m_1+\cdots+m_i+k_{i+1}+\cdots+k_t}+(-1)^{k_1-m_1+\cdots+k_i-m_i+k_{t+1}+\cdots+k_l}\right\}\frac{(m_1-1)\cdots(m_i-1)}{m_1!\cdots m_i!}\\
&\qquad\times\frac{(k_1-m_1-1)\cdots(k_i-m_i-1)(k_{i+1}-1)\cdots(k_{l}-1)}{(k_1-m_1)!\cdots(k_i-m_i)!k_{i+1}!\cdots k_{l}!}\\
&=(-1)^{l+i}\left\{(-1)^{k_{i+1}+\cdots+k_t}+(-1)^{k_1+\cdots+k_{i}+k_{t+1}\cdots+k_l}\right\}\\
&\qquad\times\prod_{j=1}^{i}\left\{\sum_{m_j=2}^{k_j-2}(-1)^{m_j}\binom{k_j}{m_j}(m_j-1)(k_j-m_j-1)\right\}\frac{(k_{i+1}-1)\cdots(k_l-1)}{k_1!\cdots k_l!}\\
&=(-1)^{l+i}\left\{(-1)^{k_{i+1}+\cdots+k_t}+(-1)^{k_1+\cdots+k_{i}+k_{t+1}\cdots+k_l}\right\}\prod_{j=1}^{i}\{1+(-1)^{k_j}\}\prod_{d=1}^{l}\frac{k_d-1}{k_d!}.
\end{align*}
Here, in the last equality, we have used the identity
 $$
 \sum_{m=2}^{k-2}(-1)^{m}\binom{k}{m}(m-1)(k-m-1)=(1+(-1)^k)(k-1).
 $$
We sum this over possible $t$ to obtain
\begin{align*}
&(-1)^{l+i}\{1+(-1)^{k_{i+1}}+\cdots+(-1)^{k_{l-1}}+(-1)^{k_l}\\
&+(-1)^{k_{i+1}+k_{i+2}}+\cdots+(-1)^{k_{l-1}+k_l}\\
&+\cdots+(-1)^{2n-k_{i+1}-k_{i+2}}+\cdots+(-1)^{2n-k_{l-1}-k_l}+\\
&\cdots+(-1)^{2n-k_{i+1}}+\cdots+(-1)^{2n-k_l}+(-1)^{2n}\}\\
&\times\prod_{j=1}^{i}\{1+(-1)^{k_j}\}\prod_{d=1}^{l}\frac{k_{d}-1}{k_{d}!}\\
&=(-1)^{l+i}\prod_{j=i+1}^{l}\left\{1+(-1)^{k_{j}}\right\}(1+(-1)^{2n})\prod_{j=1}^{i}\left\{1+(-1)^{k_j}\right\}\prod_{j=1}^{l}\frac{k_{j}-1}{k_{j}!}\\
&=(-1)^{l+i}2\prod_{j=1}^{l}\left\{1+(-1)^{k_j}\right\}\frac{k_j-1}{k_j!}.
\end{align*}

Finally, we sum this over $i$, and obtain
$$
\sum_{i=0}^{h}(-1)^{l+i}\binom{h}{i}2\prod_{j=1}^{l}\left\{1+(-1)^{k_j}\right\}\frac{k_j-1}{k_j!}=0
$$
if $h>0$ by the binomial theorem.\\

Next suppose $h=0$. Then each $k_i$ in $[k_1,\ldots,k_l]$ is either $2$ or $3$. Let $h'=\#\{i|h_i=3\}$.
The coefficient of $[k_1,\ldots,k_l]$ on the left-hand side of (\ref{eq2}) can be computed as
\begin{align*}
&(-1)^l\sum_{i=0}^{h'}\sum_{j=0}^{l-h'}(-1)^{2j+3i}\binom{h'}{i}\binom{l-h'}{j}\left(\frac{2}{3!}\right)^{h'}\left(\frac{1}{2!}\right)^{l-h'}\\
&=
\begin{dcases}
0 &\; h'\neq 0,\\
(-1)^l=(-1)^n &\; h'=0.
\end{dcases}
\end{align*}
By the binomial theorem, we have
$$
\sum_{k=0}^{2n}(-1)^kS(k)\ast S(2n-k)=(-1)^n[\underbrace{2,\cdots,2}_n].
$$
\section*{Acknowledgements}
The authors wish to thank Prof.~Masanobu Kaneko for his helpful advice. This work is supported by WISE program (MEXT) at Kyushu University.


\end{document}